\documentclass[11pt,a4paper]{article}

\usepackage[left=3.0cm,right=3.0cm,top=2.6cm,bottom=2.6cm]{geometry}

\usepackage[utf8]{inputenc}
\usepackage{microtype}
\usepackage{amsmath,amssymb,amsthm,mathtools}
\usepackage{bm}
\usepackage{algorithm}
\usepackage{algpseudocode}
\usepackage{booktabs}
\usepackage{multirow}
\usepackage{graphicx,graphics}
\graphicspath{{figures/}}
\usepackage{xcolor}
\usepackage{etoolbox}
\usepackage{enumitem}
\usepackage{caption}
\usepackage{titlefoot}
\usepackage[subnum]{cases}
\usepackage{hyperref}

\newcommand{\RR}{\mathbb{R}}

\newcommand{\vx}{\mathbf{x}}
\newcommand{\vy}{\mathbf{y}}
\newcommand{\vv}{\mathbf{v}}

\newcommand{\vg}{\mathbf{g}}

\newcommand{\vb}{\mathbf{b}}

\newcommand{\vp}{\mathbf{p}}
\newcommand{\vu}{\mathbf{u}}
\newcommand{\vw}{\mathbf{w}}

\newcommand{\vone}{\mathbf{1}}
\newcommand{\vzero}{\mathbf{0}}
\newcommand{\cG}{\mathcal{G}}
\newcommand{\cV}{\mathcal{V}}
\newcommand{\cE}{\mathcal{E}}
\newcommand{\cN}{\mathcal{N}}

\newcommand{\grad}{\nabla}
\newcommand{\norm}[1]{\left\lVert #1 \right\rVert}

\newtheorem{theorem}{Theorem}[section]

\newtheorem{assumption}{Assumption}

\numberwithin{equation}{section}

\definecolor{red}{rgb}{0.9,0,0}
\definecolor{green}{rgb}{0,0.9,0}
\definecolor{lightgreen}{rgb}{0.1,0.5,0.1}

\hypersetup{
  colorlinks=true,
  linkcolor=lightgreen,
  citecolor=lightgreen
}

\begin{document}

\title{D-ripALM: A Tuning-friendly Decentralized Relative-Type Inexact Proximal Augmented Lagrangian Method\unmarkedfntext{\hspace{-1mm}The first two authors contributed equally.}}

\author{Jiayi Zhu\thanks{School of Computer Science and Engineering, Sun Yat-sen University ({\tt zhujy86@mail2.sysu.edu.cn}).}, ~~
Hong Wang\thanks{School of Artificial Intelligence, Shenzhen Technology University ({\tt hitwanghong@163.com}) %3002 Lantian Road, Pingshan District, Shenzhen, Guangdong, China, 518118.
}, ~~
Ling Liang\thanks{Department of Mathematics, The University of Tennessee, Knoxville ({\tt liang.ling@u.nus.edu}).}, ~~
Lei Yang\thanks{(Corresponding author) School of Computer Science and Engineering, and Guangdong Province Key Laboratory of Computational Science, Sun Yat-sen University ({\tt yanglei39@mail.sysu.edu.cn}).}
}

\date{\today}

\maketitle

\begin{abstract}
This paper proposes \textbf{D-ripALM}, a \textbf{D}ecentralized \textbf{r}elative-type \textbf{i}nexact \textbf{p}roximal \textbf{A}ugmented \textbf{L}agrangian \textbf{M}ethod for consensus convex optimization over multi-agent networks. D-ripALM adopts a double-loop distributed optimization framework that accommodates a wide range of inner solvers, enabling efficient treatment of both smooth and nonsmooth objectives. D-ripALM employs a \textit{relative-type} error criterion to regulate the switching between inner and outer iterations, replacing a prescribed absolute tolerance sequence with a single relative tolerance parameter. Moreover, we establish rigorous convergence guarantees for D-ripALM under general convexity assumptions, without requiring smoothness or strong convexity conditions commonly imposed in the distributed optimization literature. Numerical experiments on synthetic and real datasets demonstrate that a restarted variant of D-ripALM achieves competitive computational efficiency and solution accuracy relative to representative decentralized algorithms.

\vspace{5mm}

\end{abstract}

%%%%%%%%%%%%%%%%%%%%%%%%%%%%%%%%%%%%%%%%%%%%%%%%%%%%%%%%%
\section{Introduction}

Driven by the rapid growth of large-scale computing systems, edge devices, and networked platforms, distributed optimization has emerged as a core paradigm in modern data-driven applications \cite{attiya2004distributed}. It is central to multi-agent sensor networks, decentralized machine learning, and smart power systems, where data are inherently distributed and centralized processing is often impractical because of communication, privacy, or robustness constraints.

In this paper, we study decentralized consensus convex optimization over a multi-agent network. Specifically, we consider a network of $n$ agents without a central coordinator. The network is represented by a communication graph $\cG = \{\cV, \cE\}$, where $\cV = \{1, \ldots, n\}$ denotes the set of agents and $\cE\subset\cV\times\cV$ denotes the set of communication links. Each agent $i \in \cV$ holds a local convex objective function $f_i: \RR^d \to (-\infty,+\infty]$, and the agents collaboratively solve the following global optimization problem
\begin{equation}\label{eq:problem}
\min_{\bm{x} \in \RR^d} \quad \sum_{i=1}^n f_i(\bm{x}),
\end{equation}
over a shared decision variable $\bm{x}\in\RR^{d}$. Since there is no central coordinator, the agents solve \eqref{eq:problem} through local computations and communication with their neighbors. In practice, each agent maintains a local copy of the decision variable and updates it using its local objective and information received from its neighbors. The agents aim to achieve consensus among these copies at an optimal solution of the global problem. Here, a communication round consists of a synchronous exchange of messages between neighbors. Summing local quantities over all agents may require several such rounds, depending on the network topology and the aggregation protocol.

Problem \eqref{eq:problem} covers a broad class of contemporary machine learning and large-scale optimization tasks. For example, in decentralized federated learning, models are trained across edge devices or data centers without sharing raw data with a central server \cite{gabrielli2023survey,ye2022decentralized,yuan2024decentralized}. In wireless sensor networks, decentralized optimization enables the fusion of distributed measurements for signal estimation and processing \cite{boyd2011distributed,ji2023distributed}. In power systems, the coordination of generators and storage units under network and operational constraints can likewise be formulated within this framework \cite{ananduta2021dalm,zhao2017consensus}. These applications highlight the pivotal role of decentralized optimization in modern large-scale computing systems.

Designing efficient distributed algorithms for \eqref{eq:problem} has been an active research topic for over two decades. A standard modeling step introduces a local copy $\vx_i \in \RR^d$ of the decision variable $\bm{x}$ for each agent $i\in\cV$ and reformulates \eqref{eq:problem} as a consensus-constrained optimization problem, namely,
\begin{equation}\label{eq:consensus-problem0}
\min_{\vx \in \RR^{nd}} \quad F(\vx) \coloneqq \sum_{i=1}^n f_i(\vx_i) \quad \mathrm{ s.t. } \quad \vx_1 = \vx_2 = \cdots = \vx_n, \\
\end{equation}
where $\vx = (\vx_1; \vx_2; \ldots; \vx_n) \in \RR^{nd}$ stacks all local copies. Building upon this consensus reformulation, researchers have developed numerous distributed algorithms by extending classical constrained optimization methods to decentralized settings. Among these approaches, penalty-based methods constitute a mainstream class. For example, the distributed gradient descent (DGD) method \cite{nedic2009distributed,yuan2016convergence} and its variants can be interpreted as gradient-type methods applied to an $\ell_2$-penalized formulation that enforces consensus among local variables. Although DGD is attractive due to its simplicity and low per-iteration cost, its convergence behavior is limited: with a fixed step size, it converges only to a neighborhood of the optimal solution; convergence to an optimal solution requires diminishing step sizes and is usually slow. Recent advances have sought to overcome these drawbacks by adopting $\ell_1$-norm penalties to enforce consensus at an optimal solution of the penalized problem \cite{peng2021byzantine,xu2018robust}, or by employing stronger subroutines to efficiently handle $\ell_2$-penalty problems with large penalty parameters \cite{zhang2021penalty}.

Another widely adopted paradigm for distributed optimization builds on the augmented Lagrangian method (ALM) \cite{hestenes1969multiplier,rockafellar1976augmented}, which augments penalty terms with Lagrange multipliers associated with the consensus constraints. Under suitable assumptions, distributed ALM-type algorithms can drive consensus and optimality residuals to zero with a fixed penalty parameter or step size. Early developments of distributed ALM-type algorithms were largely motivated by applications in signal processing and wireless communications \cite{jakovetic2011cooperative,mateos2010distributed,mota2011distributed}. These works typically adopt one of two main strategies for handling the ALM subproblems.

One strategy performs only a single (proximal) gradient step or an alternating minimization step for each subproblem, followed immediately by an update of the multiplier (dual) variable. This approach gives rise to many well-known \textit{single-loop} distributed algorithms. Representative examples include EXTRA \cite{li2020revisiting,shi2015extra}, PG-EXTRA \cite{shi2015proximal}, gradient tracking methods \cite{nedic2017achieving,nedic2023ab}, harnessing schemes \cite{qu2020harnessing}, and distributed variants of ADMM \cite{chang2014multi,shi2014linear}. While computationally lightweight, such methods often rely on restrictive smoothness assumptions and careful step-size selection to ensure stability and convergence.

The other strategy follows the classical inexact ALM framework \cite{rockafellar1976augmented}. In this approach, multiple inner iterations are carried out to approximately solve each ALM subproblem to a prescribed accuracy before updating the multiplier. This nested inner-outer structure yields a class of \textit{double-loop} distributed algorithms. Early work along this line includes \cite{jakovetic2011cooperative}, which was later extended in \cite{jakovetic2014linear} by incorporating various Jacobi- and Gauss-Seidel-type subroutines for solving the ALM subproblems. More recent representative methods include FlexPD and its variants \cite{mansoori2019general,mansoori2021flexpd}, as well as IDEAL \cite{arjevani2020ideal}. Compared with their single-loop counterparts, double-loop distributed ALM frameworks offer several notable advantages. First, they provide greater flexibility in the choice of inner solvers, enabling the use of efficient, problem-tailored algorithms. Second, empirical evidence \cite{mansoori2019general,mansoori2021flexpd} suggests that performing multiple inner iterations can substantially improve overall iteration efficiency and reduce communication overhead. Finally, from a theoretical standpoint, the classical convergence theory of ALM offers powerful analytical tools for establishing rigorous convergence guarantees for such double-loop distributed algorithms.

Despite these advantages, double-loop distributed ALM-type algorithms face a fundamental challenge: how to coordinate the inner loop (which approximately solves the ALM subproblem) with the outer loop (which updates the Lagrange multiplier), so as to guarantee convergence while maintaining efficiency. Existing studies typically address this issue by fixing the number of inner iterations or by adopting an \textit{absolute-type} error criterion to terminate the inner loop; see, e.g., \cite{arjevani2020ideal,lst2020asymptotically,liang2022qppal,liang2021inexact,rockafellar1976augmented}. Here, ``\textit{absolute}'' means that the subproblem accuracy is controlled by a prescribed tolerance sequence, typically subject to a summability requirement, rather than by the progress of the current iterate. Although conceptually simple, both strategies introduce substantial tuning burdens: their practical performance depends on the iteration count or tolerance schedule, but convergence guarantees alone do not indicate how to choose them for good practical performance. Consequently, satisfactory numerical performance often requires careful, time-consuming tuning, which can erode practical efficiency and robustness. Moreover, convergence analyses of representative methods, such as FlexPD \cite{mansoori2019general,mansoori2021flexpd} and IDEAL \cite{arjevani2020ideal}, rely on smoothness and strong convexity assumptions, which limit applicability to many real-world problems that are nonsmooth or lack strong convexity. Collectively, these limitations hinder the practicality of existing double-loop distributed ALM-type algorithms.

Motivated by recent advances in the relative-type inexact proximal augmented Lagrangian method (ripALM) \cite{yang2025ripalm,zhu2024ripalm}, we develop a \emph{tuning-friendly} double-loop distributed algorithm for decentralized consensus convex optimization. Here, tuning-friendly refers specifically to inner accuracy control through \textit{one} relative tolerance parameter, while the penalty and proximal parameters and the inner solver must still be chosen. The ripALM was originally proposed as a centralized proximal ALM framework for a broad class of constrained convex, possibly nonsmooth, problems. Its key feature is a \textit{relative-type} error criterion that controls subproblem accuracy using computable quantities tied to the current iterate, rather than a predetermined absolute tolerance. Compared with conventional absolute-type criteria, this relative-type rule eliminates the need to pre-specify a summable tolerance sequence and substantially alleviates parameter-tuning burdens. Building on these advantages, we propose a double-loop distributed ALM-type algorithm based on ripALM, termed D-ripALM. The proposed method uses the relative-type criterion to adaptively coordinate inner and outer loops while relaxing the smoothness and strong convexity assumptions commonly imposed by existing double-loop distributed ALM-type algorithms (e.g., \cite{arjevani2020ideal,mansoori2021flexpd}).

The main contributions of this work are summarized as follows:
\begin{itemize}
\item \textbf{A new decentralized proximal ALM-type algorithm with relative-type inexactness.} We develop D-ripALM, a double-loop proximal augmented Lagrangian method for decentralized consensus convex optimization. By incorporating a relative-type error criterion, D-ripALM provides an adaptive and principled mechanism for coordinating inner and outer iterations, with convergence guarantees inherited from ripALM and no need to prescribe a summable absolute tolerance sequence.

\item \textbf{A flexible, subsolver-agnostic algorithmic framework for general convex objectives.} D-ripALM handles a broad class of convex consensus problems without imposing smoothness or strong convexity assumptions. Its double-loop structure enables efficient, problem-tailored inner solvers. Numerical experiments in Section \ref{sec:experiments} examine the effect of tolerance choices and compare an implementation with an auxiliary-variable restart heuristic against representative double-loop and single-loop methods.
\end{itemize}

% The main contributions of this work are summarized as follows:
% \begin{itemize}
%     \item We successfully extend the centralized ripALM framework to decentralized consensus optimization, termed as D-ripALM. The proposed method inherits both the theoretical and practical advantages of ripALM, including rigorous convergence guarantees, practical improvements in parameter tuning, and robust numerical performance.
%     \item The proposed D-ripALM provides a flexible distributed algorithmic framework that is capable of handling general convex optimization problems. Its double-loop structure permits the selection of efficient and problem-tailored subroutines (e.g., FISTA \cite{beck2009fast}) to achieve efficient solution. Numerical experiments in Section \ref{sec:experiments} demonstrate the effectiveness and robustness of the proposed approach.
% \end{itemize}

The remainder of this paper is organized as follows. Section~\ref{sec:setup} introduces the problem formulation and basic assumptions. Section~\ref{sec:d-ripalm} presents D-ripALM and its distributed implementation. Section~\ref{sec:convergence} establishes the convergence properties of D-ripALM. Section \ref{sec:experiments} evaluates the numerical performance of the restarted implementation D-ripALM-re. Finally, Section~\ref{sec:conclusion} concludes the paper and outlines directions for future research.

\paragraph{Notation.} Scalars, vectors, and matrices are denoted by lowercase letters, bold lowercase letters, and uppercase letters, respectively. We use $\mathbb{R}^n$ ($\mathbb{R}_{+}^n$, $\mathbb{R}_{++}^n$) and $\mathbb{R}^{m \times n}$ to denote the sets of $n$-dimensional (non-negative, positive) real vectors and $m \times n$ real matrices, respectively. $\vzero_{n}$ (or $\bm{1}_{n}$) denotes the all-zero (or all-one) vector in $\mathbb{R}^{n}$, $I_{n} \in \mathbb{R}^{n\times n}$ denotes the identity matrix, and $O_{m\times n}\in\mathbb{R}^{m\times n}$ denotes the zero matrix. When the dimension is clear, we omit subscripts (i.e., $\vzero$, $\bm{1}$, $I$, $O$). $\otimes$ denotes the Kronecker product. For a vector $\bm{x}\in\RR^{n}$, $x_{i}$ denotes its $i$-th entry, $\|\bm{x}\|$ denotes its Euclidean norm, $\|\bm{x}\|_{1}$ denotes its $\ell_{1}$-norm defined by $\|\bm{x}\|_{1} \coloneqq \sum_{i=1}^{n} |x_{i}|$, and $\|\bm{x}\|_{\infty}$ denotes its $\ell_{\infty}$-norm given by the largest entry in magnitude.

For an extended-real-valued function $f: \mathbb{R}^{n} \rightarrow [-\infty,\infty]$, we say $f$ is \textit{proper} if $f(\bm{x}) > -\infty$ for all $x\in\mathbb{R}^{n}$ and its effective domain ${\rm dom}(f)\coloneqq\{\bm{x} \in \mathbb{R}^{n} \mid f(\bm{x})<\infty\}$ is nonempty. A proper function $f$ is said to be closed if it is lower semicontinuous. Let $f: \mathbb{R}^{n} \rightarrow (-\infty,\infty]$ be proper, closed, and convex; the subdifferential of $f$ at $\bm{x}\in{\rm dom}(f)$ is defined by $\partial f(\bm{x})\coloneqq\big\{\bm{d}\in\mathbb{R}^n \mid f(\bm{y}) \geq f(\bm{x}) + \langle \bm{d}, \,\bm{y}-\bm{x}\rangle, ~\forall\,\bm{y}\in\mathbb{R}^n\big\}$. For any $\nu>0$, the proximal mapping of $\nu f$ at $\bm{x}$ is defined by $\mathtt{prox}_{\nu f}(\bm{x}) \coloneqq \arg\min_{\bm{y}} \big\{f(\bm{y}) + \frac{1}{2\nu}\|\bm{y} - \bm{x}\|^2\big\}$. Let $\mathcal{S}$ be a closed convex subset of $\mathbb{R}^n$. Its indicator function $\delta_{\mathcal{S}}$ is defined by $\delta_{\mathcal{S}}(\bm{x})=0$ if $x\in\mathcal{S}$ and $\delta_{\mathcal{S}}(\bm{x})=+\infty$ otherwise. Moreover, we denote the Euclidean distance of $\bm{x}\in\mathbb{R}^{n}$ to $\mathcal{S}$ by $\mathrm{dist}(\bm{x},\mathcal{S})\coloneqq\inf_{\bm{y}\in\mathcal{S}}\|\bm{x} - \bm{y}\|$.

%%%%%%%%%%%%%%%%%%%%%%%%%%%%%%%%%%%%%%%%%%%%%%%%%%%%%%%%%%%%%%%%%%%%%%%%%
\section{Problem setup}\label{sec:setup}

%%%%%%%%%%%%%%%%%%%%%%%%%%%%%%%%%%%%%%%%%%%%%%%%%%%%%%%%%%%%%%%%%%%%%%%%%
\subsection{Assumptions on the objective functions and network}

We make the following blanket assumptions on problem \eqref{eq:problem}.

\begin{assumption}\label{ass:smoothness}
(i) Each function $f_i:\RR^d\to(-\infty,+\infty]$ in \eqref{eq:problem} is proper, closed, and convex (possibly nonsmooth); (ii) $\bigcap_{i=1}^{n}\mathrm{dom}(f_i)\neq\emptyset$; (iii) problem \eqref{eq:problem} admits at least one optimal solution $\bm{x}^*$.
\end{assumption}

Under Assumption \ref{ass:smoothness}, the aggregate objective $f(\bm{x})\coloneqq\sum_{i=1}^n f_i(\bm{x})$ is proper, closed, and convex with a nonempty effective domain. Neither differentiability nor Lipschitz continuity of $\grad f_i$ is required, and no strong convexity of $f$ is assumed.
%Here, when $f_i$ is smooth, we denote its gradient by $\grad f_i$; otherwise, we select a certain subgradient $\vg_i(\vx_i)\in\partial f_i(\vx_i)$.
These assumptions are consistent with those adopted in the centralized ripALM framework \cite{yang2025ripalm,zhu2024ripalm} and enable a unified treatment of both smooth and nonsmooth objectives. In practice, solving the resulting subproblems requires access to a first-order oracle for each $f_i$, either via gradient evaluations or proximal mappings.

\begin{assumption}\label{ass:mixing-mat}
Consider a connected communication graph $\cG=\{\mathcal{V},\mathcal{E}\}$ consisting of a set of agents $\cV=\{1,\ldots,n\}$ and a set of undirected edges $\cE\subset \mathcal{V}\times \mathcal{V}$. The associated mixing matrix $W=[w_{ij}] \in \mathbb{R}_{+}^{n \times n}$ satisfies: (i) decentralization: $w_{ij} = 0$ if $(i,j) \notin \cE$ and $i \neq j$; (ii) symmetry: $W = W^{\mathsf{T}}$; (iii) double stochasticity: $W \vone_n = \vone_n$ and $\vone_n^{\mathsf{T}} W = \vone_n^{\mathsf{T}}$; (iv) spectral bounds: $-I_n \prec W \preceq I_n$; (v) consensus null space: $\ker(I_{n}-W)=\operatorname{span}\{\vone_{n}\}$.
% \begin{enumerate}[label=(\roman*)]
% % \item \emph{Connectivity}: $\cG$ is connected and undirected;
% \item \emph{Decentralized property}: $w_{ij} = 0$ if $(i,j) \notin \cE$ and $i \neq j$.
% \item \emph{Symmetry}: $W = W^{\mathsf{T}}$.
% \item \emph{Doubly stochastic}: $W \vone_n = \vone_n$ and $\vone_n^{\mathsf{T}} W = \vone_n^{\mathsf{T}}$.
% \item \emph{Spectral bounds}: $-I_n \prec W \preceq I_n$.
% \end{enumerate}
\end{assumption}

Assumption \ref{ass:mixing-mat} guarantees that the mixing matrix $W$ has real eigenvalues satisfying $1 = \lambda_1(W) > \lambda_2(W) \ge \cdots \ge \lambda_n(W) > -1$. The spectral gap, defined as $1-\varsigma$ with $\varsigma \coloneqq \max\{|\lambda_2(W)|, |\lambda_n(W)|\} < 1$, quantifies the connectivity of the underlying communication network. The associated network condition number $\kappa_W\coloneqq \frac{1+\varsigma}{1-\varsigma}$ characterizes the conditioning of the consensus process. Smaller values of $\varsigma$ (equivalently, a larger spectral gap) correspond to better-connected networks and faster consensus dynamics.

%%%%%%%%%%%%%%%%%%%%%%%%%%%%%%%%%%%%%
\subsection{Consensus reformulation}

Under Assumption \ref{ass:mixing-mat}, the consensus constraint $\vx_1=\cdots=\vx_n$ in \eqref{eq:consensus-problem0} is equivalent to the linear constraint
\begin{equation*}
Z\vx = \vzero, \quad \text{ where }\quad Z \coloneqq (I_n - W) \otimes I_d \in \RR^{nd \times nd},
\end{equation*}
where $W$ is the mixing matrix associated with the communication graph.
%Several basic properties of $Z$ follow directly from Assumption \ref{ass:mixing-mat}.
Because $I_n - W \succeq 0$ and the Kronecker product preserves positive semidefiniteness, $Z$ is positive semidefinite. Moreover, its null space coincides with the consensus subspace, i.e.,
\begin{equation}\label{eq:ker-Z}
\mathrm{null}(Z) = \mathrm{span}\{\vone_n \otimes \bm{v} \mid \bm{v} \in \RR^d\},
\end{equation}
which implies that $Z\vx = \vzero$ if and only if $\vx_1 = \cdots = \vx_n$. The spectral properties of $Z$ are inherited from those of $W$. Let $W = Q \Lambda_W Q^\mathsf{T}$ be an eigenvalue decomposition of $W$, then
\begin{equation*}
Z = (Q \otimes I_d)((I_n - \Lambda_W) \otimes I_d)(Q^\mathsf{T} \otimes I_d),
\end{equation*}
so the eigenvalues of $Z$ are $\{1-\lambda_i(W) \mid i=1,\ldots,n\}$, each with multiplicity $d$. In particular, the smallest and largest positive eigenvalues of $Z$ are
$\lambda_{\min}^+(Z)=1-\lambda_2(W)$ and $\lambda_{\max}(Z)=1-\lambda_n(W)$, which govern the consensus rate.

Since $Z$ is positive semidefinite, there exists a symmetric square root $\sqrt{Z}$ such that $Z = (\sqrt{Z})^{2}$. Because $\sqrt{Z}$ shares the same null space as $Z$, the constraint $Z\vx=\vzero$ is equivalent to $\sqrt{Z}\vx = \vzero$. Using this relation, problem \eqref{eq:consensus-problem0} can be reformulated as
\begin{equation}\label{eq:consensus-problem}
\min_{\vx \in \RR^{nd}}\quad F(\vx) \coloneqq \sum_{i=1}^n f_i(\vx_i) \quad \text{s.t.}\quad \sqrt{Z}\vx = \vzero.
\end{equation}

%%%%%%%%%%%%%%%%%%%%%%%%%%%%%%%%%%%%%%%%%%%%%%%%%%%%%%%%%%%%%%%%%%%%%%%%%%%%%%%
\section{A decentralized relative-type inexact proximal ALM}\label{sec:d-ripalm}

We now develop a decentralized relative-type inexact proximal augmented Lagrangian method for solving the consensus problem \eqref{eq:consensus-problem}. Relative-type error criteria have been well studied in centralized settings \cite{yang2025ripalm,zhu2024ripalm}. Here, we adapt these ideas to decentralized environments while ensuring that all algorithmic steps can be implemented in a fully distributed manner.

Given a penalty parameter $\sigma>0$, the augmented Lagrangian associated with \eqref{eq:consensus-problem} is
\begin{equation}\label{eq:auglag-natural}
\mathcal{L}_\sigma(\vx, \vy) \coloneqq F(\vx) + \langle \vy, \sqrt{Z}\vx \rangle + \frac{\sigma}{2}\norm{\sqrt{Z}\vx}^2,
\end{equation}
where $\vy = (\vy_{1};\vy_{2};\ldots;\vy_{n} ) \in \RR^{nd}$ denotes the Lagrange multiplier associated with the constraint $\sqrt{Z}\vx=\vzero$. Starting from an initial point $(\vx^{0},\vy^{0}) \in \mathrm{dom}(F)\times\RR^{nd}$, the classical inexact ALM for \eqref{eq:consensus-problem} is
\begin{numcases}{}
\mathbf{x}^{k+1} \approx \arg\min\limits_{\mathbf{x}\in\mathbb{R}^{nd}} ~ \left\{\mathcal{L}_{\sigma}(\mathbf{x}, \vy^{k}) \right\}, \label{eq:ALM-subpro}\\
\vy^{k+1} = \vy^{k} + \sigma \sqrt{Z}\vx^{k+1}, \label{eq:ALM-multi}
\end{numcases}
which alternates between approximately minimizing the ALM subproblem \eqref{eq:ALM-subpro} and updating the multiplier in \eqref{eq:ALM-multi}. Following the centralized ripALM framework developed in \cite{yang2025ripalm,zhu2024ripalm}, we enhance the classical ALM with two key ingredients: (i) an additional proximal term $\frac{\tau_{k}}{2\sigma_{k}}\|\vx-\vx^{k}\|^{2}$, which enforces strong convexity of the subproblem and improves its conditioning; and (ii) a relative-type error criterion that adaptively regulates the accuracy required for solving each subproblem as the algorithm progresses. The resulting decentralized algorithm for \eqref{eq:consensus-problem}, termed D-ripALM, is summarized in Algorithm~\ref{alg:D-ripALM}.

\begin{algorithm}[ht]
\caption{D-ripALM for solving problem \eqref{eq:consensus-problem}}\label{alg:D-ripALM}
\textbf{Input:} $\rho\in[0,1)$, $\{\sigma_{k}\}_{k=0}^{\infty}\subset\mathbb{R}_{++}$ with $\inf_{k\geq0}\{\sigma_{k}\}>0$, proximal parameter sequence $\{\tau_k\}_{k=0}^\infty \subset \RR_{++}$ with $\inf_{k \ge 0} \tau_k > 0$ and $\sup_{k \ge 0} \tau_k < \infty$. Choose $\vw^0\in\mathbb{R}^{nd}$ arbitrarily, choose $\vx^0\in\mathrm{dom}(F)$ and let $\vy^{0}=\vzero$. Set $k=0$.  \vspace{1mm} \\
\textbf{while} the termination criterion is not met, \textbf{do} \vspace{-2mm}
	\begin{itemize}[leftmargin=2.2cm]
	\item[\textbf{Step 1}.] Starting from $\vx^{k}$, apply a suitable decentralized method to solve
		\begin{equation*}
		\min\limits_{\vx\in\mathbb{R}^{nd}}~~
		\mathcal{L}_{\sigma_{k}}(\vx,\,\vy^k)
		+ \frac{\tau_k}{2\sigma_{k}}\big\|\vx-\vx^k\big\|^2,
		\end{equation*}
		to find an approximate solution $\vx^{k+1}$ and an error term $\Delta^{k+1}$
		such that
		\begin{equation}\label{eq:D-ripALM-errterm}
		\Delta^{k+1}\in
		\partial_{x}\mathcal{L}_{\sigma_{k}}(\vx^{k+1},\,\vy^k)
            + \tau_k\sigma_{k}^{-1}\big(\vx^{k+1}-\vx^{k}\big),
		\end{equation}
		satisfying the following relative-type error criterion
		\begin{equation}\label{eq:D-ripALM-errcrit}
		\begin{aligned}
		&{}2\big|\langle\vw^k-\vx^{k+1},
		\,\sigma_{k}\Delta^{k+1}\rangle\big|
		+ \big\|\sigma_{k}\Delta^{k+1}\big\|^2 \\
		\leq{}& \rho\left(
		\big\|\sigma_{k}\sqrt{Z}\vx^{k+1}\big\|^2
		+ \tau_k\big\|\vx^{k+1}-\vx^k\big\|^2\right).
    \end{aligned}
		\end{equation}
		
		\item[\textbf{Step 2}.] Update
		\begin{equation}\label{eq:D-ripALM-ywupdate}
		\vy^{k+1} = \vy^k + \sigma_{k}\sqrt{Z}\vx^{k+1},
            \qquad
            \vw^{k+1} = \vw^k - \sigma_{k}\Delta^{k+1}.
		\end{equation}
		
		\item[\textbf{Step 3}.] Set $k=k+1$ and go to \textbf{Step 1}.  \vspace{-2mm}
	\end{itemize}
\textbf{end while}  \\
\textbf{Output:} $(\vx^k,\vy^k)$ \vspace{0.5mm}
\end{algorithm}

We now describe how Algorithm~\ref{alg:D-ripALM} can be implemented in a fully distributed manner. At the $k$-th outer iteration, the agents cooperatively solve the subproblem
\begin{equation}\label{eq:D-ripALM-subpro}
\min_{\vx\in\RR^{nd}} \quad
\Psi_{k}(\vx) \coloneqq F(\vx)
+ \langle \Omega^{k},\vx\rangle
+ \frac{\sigma_{k}}{2}\langle\vx,Z\vx\rangle
+ \frac{\tau_{k}}{2\sigma_{k}}\big\|\vx-\vx^{k}\big\|^{2},
\end{equation}
where $\Omega^{k}=\sqrt{Z}\vy^{k}$ is a transformed dual variable. This transformation avoids the explicit appearance of $\sqrt{Z}$, which is not directly accessible in decentralized settings, while still implicitly tracking the original dual variable $\vy^{k}$. The proximal term $\frac{\tau_{k}}{2\sigma_{k}}\|\vx-\vx^{k}\|^{2}$ makes $\Psi_k$ strongly convex and thus guarantees a unique minimizer of \eqref{eq:D-ripALM-subpro} even when $\mathcal{L}_{\sigma}(\cdot, \vy^{k})$ is merely convex. Moreover, this proximal term improves the conditioning of the subproblem and facilitates efficient inner minimization using problem-tailored subsolvers.% under the proposed relative-type error criterion \eqref{eq:D-ripALM-errcrit}.

Rather than solving \eqref{eq:D-ripALM-subpro} exactly, D-ripALM accepts an inexact solution $\vx^{k+1}$ that satisfies the relative-type error criterion \eqref{eq:D-ripALM-errcrit}. In practice, (proximal) gradient methods and their accelerated variants are well suited as inner solvers because of their computational efficiency, compatibility with nonsmooth objectives, and amenability to distributed implementations. Specifically, gradient-type methods require repeated evaluations of the gradient of the smooth part of $\Psi_k$, which includes the quadratic penalty term
\begin{equation*}
\frac{\sigma_k}{2}\nabla_\vx\langle\vx, Z\vx\rangle = \frac{\sigma_k}{2}(Z + Z^\mathsf{T})\vx = \sigma_k Z\vx.
\end{equation*}
This term is locally computable after a single round of neighbor communication. Indeed, each agent $i$ can form $[Z\vx]_i = \vx_i - \sum_{j\in\cN_i\cup\{i\}}w_{ij}\vx_j$ using only information received from its neighbors. Consequently, gradient-type methods can be implemented efficiently in a fully distributed manner through local computations and neighbor communications.

We next explain how to verify the relative-type error criterion \eqref{eq:D-ripALM-errcrit} in a distributed manner. Because the variables in the distributed problem are block structured, the quantities $\vw^{k}$, $\vx^{k}$, $\vx^{k+1}$, $\Delta^{k+1}$, $\Omega^{k}$, and $Z\vx^{k+1}$ are block-wise separable, and the $i$-th block of each is held by agent $i$. Using this structure, the error criterion \eqref{eq:D-ripALM-errcrit} can be rewritten as the following finite-sum inequality across agents:
\begin{equation}\label{eq:D-ripALM-errterm2}
	\begin{aligned}
		& 2\left|\sum_{i=1}^{n}\langle\vw_{i}^k-\vx_{i}^{k+1},
		\,\sigma_{k}\Delta^{k+1}_{i}\rangle\right|
		+ \sum_{i=1}^{n}\big\|\sigma_{k}\Delta^{k+1}_{i}\big\|^2 \\
		  \leq  &\;\; \rho\sum_{i=1}^{n}\left(
		\sigma_{k}^{2}\langle \vx_{i}^{k+1}, [Z\vx^{k+1}]_{i} \rangle
		+ \tau_k\big\|\vx_{i}^{k+1}-\vx_{i}^k\big\|^2\right).
	\end{aligned}
\end{equation}
Thus, verifying \eqref{eq:D-ripALM-errcrit} reduces to checking whether \eqref{eq:D-ripALM-errterm2} holds, which requires collaboration among agents. Using the summation form and block separability, a distributed verification strategy for \eqref{eq:D-ripALM-errterm2} proceeds as follows: at the $k$-th outer iteration, each agent computes the following local variable stack
\begin{align*}
    E_{i}(k) \coloneqq \left[E_{i}^{1}(k), E_{i}^{2}(k), E_{i}^{3}(k)\right],
\end{align*}
where $E_{i}^{1}(k) = \langle\vw_{i}^{k} - \vx_{i}^{k+1}, \sigma_k\Delta_{i}^{k+1}\rangle$, $E_{i}^{2}(k) = \|\sigma_k\Delta_{i}^{k+1}\|^{2}$ and
\begin{equation*}
E_{i}^{3}(k) = \sigma_{k}^{2}\langle\vx_{i}^{k+1}, [Z\vx^{k+1}]_{i}\rangle + \tau_{k}\|\vx_{i}^{k+1} - \vx_{i}^{k}\|^2.
\end{equation*}
In particular, when gradient-type methods are used as subsolvers, the $i$-th block of the error term $\Delta^{k+1}$ can be computed as
\begin{equation}\label{eq:grad-natural}
\Delta_{i}^{k+1} \coloneqq \vg_{i}(\vx_{i}^{k+1}) + \Omega^{k}_{i} + \sigma_{k} [Z \vx^{k+1}]_{i} + \tau_k\sigma_k^{-1}\big(\vx_{i}^{k+1} - \vx_{i}^k\big),
\end{equation}
where $\vg_{i}(\vx_{i}^{k+1}) = \grad f_{i}(\vx_{i}^{k+1})$ if $f_{i}$ is differentiable, and $\vg_{i}(\vx_{i}^{k+1})\in\partial f_{i}(\vx_{i}^{k+1})$ otherwise.

Once the local quantities $E_{i}(k)$ have been computed, the agents aggregate them through neighbor communication and make the resulting three sums available to all agents. The candidate $\vx^{k+1}$ is accepted if
\begin{equation*}
2\left|\sum_{i=1}^{n}E_{i}^{1}(k)\right| + \sum_{i=1}^{n}E_{i}^{2}(k) \leq \rho \sum_{i=1}^{n}E_{i}^{3}(k).
\end{equation*}
Each aggregation message contains only three scalars, representing local contributions or accumulated sums. For large $d$, each such message therefore contains far fewer scalars than a $d$-dimensional vector exchanged by the inner solver. However, completing the check may require multiple communication rounds to pass information through intermediate agents and make the sums available throughout the network. The resulting communication and synchronization costs depend on the network topology and the aggregation protocol.

Finally, we note that the update of the ordinary dual variable $\vy^{k+1}$ in \eqref{eq:D-ripALM-ywupdate} involves the square-root matrix $\sqrt{Z}$, which is not directly accessible in decentralized settings. However, the algorithm only requires $\Omega^{k+1}=\sqrt{Z}\vy^{k+1}$ for all $k\geq0$, and this quantity is updated by
\begin{equation}\label{eq:dual-relation}
	\Omega^{k+1} = \Omega^k + \sigma_k Z \vx^{k+1}.
\end{equation}
Therefore, $\vy^{k+1}$ does not need to be explicitly computed. Moreover, since $Z\vx^{k+1}$ has already been computed and cached during the verification of \eqref{eq:D-ripALM-errcrit}, no additional communication is required for the dual update.

%%%%%%%%%%%%%%%%%%%%%%%%%%%%%%%%%%%%%%%%%%%%%%%%%%%%%%%%%%%%%%%%%%%%%%%%%
\section{Convergence analysis}\label{sec:convergence}

In this section, we study the convergence properties of D-ripALM in Algorithm \ref{alg:D-ripALM} under Assumptions \ref{ass:smoothness}--\ref{ass:mixing-mat}. %The main theoretical results are adapted from the convergence analysis to the centralized ripALM in \cite{yang2025ripalm}, since problem \eqref{eq:consensus-problem} can be viewed as a special case of the problem formulation considered in \cite{yang2025ripalm}.
For problem \eqref{eq:consensus-problem}, the Lagrangian function is defined as
\begin{equation}\label{eq:lagrangian-consensus}
\ell(\vx, \vy) \coloneqq F(\vx) + \langle \vy, \sqrt{Z}\vx \rangle,
\end{equation}
and its subdifferential is given by
\begin{equation*}\label{eq:subdiff-lag}
\partial \ell(\vx, \vy) = \big\{\partial F(\vx) + \sqrt{Z}\vy\big\} \times \big\{ - \sqrt{Z}\vx \big\}.
\end{equation*}
We call $(\vx^*, \vy^*)$ a \textit{saddle point} of $\ell$ if $(\vzero, \vzero) \in \partial \ell(\vx^*, \vy^*)$, i.e.,
\begin{numcases}{}
\vzero = \vv^{*} + \sqrt{Z}\vy^*, \quad  \exists\,\vv^{*}\in \partial F(\vx^*), \label{eq:stationary}\\[3pt]
-\sqrt{Z}\vx^* = \vzero. \label{eq:consensus}
\end{numcases}

With the above preparations, we are ready to state the convergence results for D-ripALM in Algorithm \ref{alg:D-ripALM}. Notably, D-ripALM can be viewed as a distributed implementation of \cite[Algorithm 1]{yang2025ripalm}, since the consensus problem \eqref{eq:consensus-problem} is a linearly constrained convex optimization problem that fits the template studied in \cite{yang2025ripalm}. Therefore, all theoretical results in \cite[Section 3]{yang2025ripalm} apply directly. For brevity, we omit the proofs here and refer readers to \cite{yang2025ripalm} for detailed derivations.

\begin{theorem}\label{thm:global-convergence}
Let $\rho\in[0,1)$, $\{\sigma_{k}\}$ be a positive sequence satisfying $\sigma_k\geq\sigma_{\min}>0$ for all $k\geq0$, and $\{\tau_k\}$ be a positive sequence satisfying
\begin{equation*}
\tau_{k}\geq\tau_{\min}>0, \quad \tau_{k+1}\leq(1+\nu_{k})\tau_{k} \quad
\mbox{with} \quad \nu_{k}\geq0 ~~\mbox{and}~~ {\textstyle\sum_{k=0}^{\infty}}\nu_{k} < +\infty.
\end{equation*}
Let $\{\vx^{k}\}$, $\{\Delta^{k}\}$, $\{\vw^{k}\}$, and $\{\vy^{k}\}$ be the iterates generated by D-ripALM in Algorithm \ref{alg:D-ripALM}. If $\ell$ admits a saddle point, i.e., $(\partial\ell)^{-1}(\bm{0},\bm{0})\neq\emptyset$, then:

\begin{enumerate}[label=(\roman*),leftmargin=*]
\item The sequences $\{\vx^k\}$, $\{\vw^k\}$, and $\{\vy^k\}$ are bounded.

\item $\lim_{k \to \infty} \Delta^{k+1} = \vzero$, $\lim_{k \to \infty} \vp^{k+1} = \vzero$, and $\lim_{k \to \infty} \vu^{k+1} = \vzero$, where
\begin{equation*}
\vp^{k+1} \coloneqq \Delta^{k+1} - \tau_k\sigma_k^{-1}(\vx^{k+1} - \vx^k) \quad \text{and} \quad \vu^{k+1} \coloneqq \sigma_k^{-1}(\vy^k - \vy^{k+1}),\quad \forall k\geq 0.
\end{equation*}

\item The sequence $\{F(\vx^{k+1})\}$ converges to the optimal value of \eqref{eq:consensus-problem}.

\item Any accumulation point of $\{\vx^k\}$ is an optimal solution of problem \eqref{eq:consensus-problem}, and any accumulation point of $\{\vy^k\}$ is a corresponding optimal dual multiplier.

\item The entire sequence $\{\vy^k\}$ converges to an optimal dual multiplier $\vy^*$.
\end{enumerate}
\end{theorem}

We next present the asymptotic (super)linear convergence rate for the outer iterates of D-ripALM under an additional error bound condition.

\begin{assumption}\label{ass:error-bound}
For any $r > 0$, there exists a constant $\kappa > 0$ such that for all $(\vx, \vy) \in \RR^{nd} \times \RR^{nd}$ with $\mathrm{dist}((\vx, \vy), (\partial \ell)^{-1}(\vzero, \vzero)) \le r$,
\begin{equation}\label{eq:error-bound}
\mathrm{dist}((\vx, \vy), (\partial \ell)^{-1}(\vzero, \vzero)) \le \kappa \cdot \mathrm{dist}((\vzero, \vzero), \partial \ell(\vx, \vy)).
\end{equation}
\end{assumption}

As discussed in \cite[Section 3.1]{yang2025ripalm}, imposing a suitable error bound condition is a standard approach for deriving fast asymptotic convergence rates of ALM-type algorithms in convex optimization \cite{lst2020asymptotically,rockafellar1976augmented,yang2025ripalm,zhu2024ripalm}. In Rockafellar's seminal works \cite{rockafellar1976augmented,r1976monotone}, this was achieved under the assumption that $(\partial\ell)^{-1}$ is Lipschitz continuous at the origin. However, this assumption is restrictive because it implicitly assumes uniqueness of the optimal solution. To relax this condition, Luque \cite{l1984asymptotic} introduced a growth condition, which is essentially equivalent to local upper Lipschitz continuity \cite{r1976implicit,r1981some}. Regarding the error bound condition \eqref{eq:error-bound} in Assumption \ref{ass:error-bound}, \cite[Lemma 2.4]{lst2020asymptotically} indicates that it is even weaker than Luque's error bound in \cite{l1984asymptotic}. We now revisit the asymptotic fast convergence rate of D-ripALM under Assumption \ref{ass:error-bound}.

\begin{theorem}[\textbf{Asymptotic (super)linear convergence}]\label{thm:superlinear}
Let $\ell$ be defined as in \eqref{eq:lagrangian-consensus}, and let $\rho\in[0,1)$, $\{\sigma_{k}\}$ be a positive sequence satisfying $\sigma_k\geq\sigma_{\min}>0$ for all $k\geq0$, and $\{\tau_{k}\}$ be a positive sequence satisfying
\begin{equation*}
\tau_{k}\geq\tau_{\min}>0,
\quad \tau_{k+1}\leq(1+\nu_{k})\tau_{k}
\quad \mbox{with} \quad \nu_{k}\geq0
~~\mbox{and}~~
{\textstyle\sum_{k=0}^{\infty}}\nu_{k} < +\infty.
\end{equation*}
Suppose additionally that $\ell$ admits a saddle point (i.e., $(\partial\ell)^{-1}(\vzero,\vzero)\neq\emptyset$), Assumption \ref{ass:error-bound} holds, and the parameters $\rho$, $\{\sigma_k\}$, and $\{\tau_{k}\}$ satisfy
\begin{equation}\label{para-conds}
\sqrt{\tau_{\min}} - 2\sqrt{\rho} > 0
\quad \text{and} \quad
\liminf\limits_{k\to\infty} ~ \sigma_{k} > c\cdot\frac{2\kappa\sqrt{\tau_{\max}}\left(\rho+\sqrt{\rho\,\overline{\tau}_{\max}}\right)}{\sqrt{\tau_{\min}} - 2\sqrt{\rho}},
\end{equation}
where $c>1$ is an arbitrarily given positive constant, $\tau_{\max}\coloneqq\tau_{0} \prod_{k=0}^{\infty}(1+\nu_{k})$, and $\overline{\tau}_{\max}\coloneqq\max\left\{1,\tau_{\max}\right\}$. Let $\Lambda_k \coloneqq \mathrm{Diag}(\tau_{k}I_{nd},I_{nd})$, $\overline{\tau}_{k} \coloneqq \max\left\{1, \tau_{k}\right\}$, and
\begin{equation*}
\gamma_{k} \coloneqq \left(1 - \frac{2\kappa\sqrt{\tau_{k}}\left(\rho+\sqrt{\rho\overline{\tau}_{k}}\right) + 2\sigma_{k}\sqrt{\rho}}{\sigma_{k}\sqrt{\tau_{k}}}\right) \frac{\sigma_{k}^{2}}{\kappa^2\left(\sqrt{\rho}+\sqrt{\overline{\tau}_{k}}\right)^2\overline{\tau}_{k}}.
\end{equation*}
Then the following statements hold.
\begin{enumerate}[label=(\roman*)]
\item For all sufficiently large $k$, we have
		\begin{equation*}
			\gamma_{k}\geq
			\left(\frac{c-1}{c}\right)
			\cdot\frac{\sqrt{\tau_{\min}} - 2\sqrt{\rho}}{\sqrt{\tau_{\min}}}
			\cdot\frac{\sigma_k^2}{\kappa^2\left(\sqrt{\rho}+\sqrt{\overline{\tau}_{\max}}\right)^2\overline{\tau}_{\max}} > 0,
		\end{equation*}
		and
		\begin{equation*}%\label{eq:linear_rate}
			\mathrm{dist}_{\Lambda_{k+1}}\left((\vx^{k+1},\vy^{k+1}),
			\,(\partial\ell)^{-1}(\vzero, \vzero)\right) \leq \mu_{k}\, \mathrm{dist}_{\Lambda_k}\left((\vx^{k},\vy^{k}), \,(\partial\ell)^{-1}(\vzero, \vzero)\right),
		\end{equation*}
		where
		\begin{equation*}
			%\limsup\limits_{k\to\infty}\,\left\{\mu_{k}\coloneqq\sqrt{\frac{1+\nu_{k}}{1+\gamma_{k}}}\right\} < 1.
			\mu_{k}\coloneqq\sqrt{\frac{1+\nu_{k}}{1+\gamma_{k}}}
			\quad \mbox{satisfies} \quad
			\limsup\limits_{k\to\infty}\,\left\{\mu_{k}\right\}<1
			~~\mbox{as}~~\nu_{k}\to0.
		\end{equation*}
		
\item The entire sequence $\{\vx^k\}$ converges.
\end{enumerate}
\end{theorem}

% Before we close this section, we remark here that the optimal solution to problem \eqref{eq:consensus-problem} induces the optimal solution to \eqref{eq:problem}. According to \eqref{eq:consensus}, we know that $\vx^{*}$ is consensus, and there must exist a vector $x^{*}\in\RR^{d}$ such that $\vx^{*} = \vone_{n}\otimesx^{*}$ according to \eqref{eq:ker-Z}. Then, multiplying both sides of \eqref{eq:stationary} by $\bm{1}_{n}^{\mathsf{T}}\otimes I_{d}$, we derive from \eqref{eq:ker-Z} that $(\vone_{n}^{\mathsf{T}}\otimes I_{d})\sqrt{Z}\vy^{*} = \vzero$, and consequently $(\vone_{n}^{\mathsf{T}}\otimes I_{d})\vv^{*} = \vzero$ as well. Finally, substituting $\vx^{*} = \vone_{n}\otimesx^{*}$ and according to that $\vv^{*}\in\partial F(\vx^{*})$, we conclude that
% \begin{equation*}
% \bm{0}\in\sum_{i=1}^{n} \partial f_{i}(x^{*}),
% \end{equation*}
% which indicates that $x^{*}$ is the optimal solution to problem \eqref{eq:problem}.

%%%%%%%%%%%%%%%%%%%%%%%%%%%%%%%%%%%%%%%%%%%%%%%%%%%%%%%%%%%%%%%%%%%%%%%%%
\section{Numerical experiments}\label{sec:experiments}

In this section, we evaluate D-ripALM-re, a variant of D-ripALM (Algorithm \ref{alg:D-ripALM}) that incorporates the auxiliary-variable restart strategy described below, for solving two distributed optimization problems of the form \eqref{eq:problem}. Specifically, the experiments are organized as follows.

\begin{itemize}
\item In Section \ref{sec:comp_diff_paras}, we compare D-ripALM-re with IDEAL \cite{arjevani2020ideal}, a double-loop distributed ALM-type algorithm using an absolute-type error criterion, for solving the $\ell_2$-regularized logistic regression problem \eqref{eq:log-reg} using synthetic data and the MNIST dataset. We examine the effects of tolerance choices on numerical performance.

\item In Section \ref{sec:comp_diff_algos}, we compare D-ripALM-re with PG-EXTRA \cite{shi2015proximal} and NIDS \cite{li2019nids}, two representative single-loop decentralized methods, for solving the $\ell_1$-regularized least-squares regression (LASSO) problem \eqref{eq:lasso} using synthetic data and real datasets from the LIBSVM collection. We examine the computational costs and accuracy attained under a common iteration budget across different regularization parameters and network topologies.

\end{itemize}

All experiments are conducted in {\sc Matlab} R2023a on a PC with Intel processor i7-12700K@3.60GHz (with 12 cores and 20 threads), 128GB of RAM and Windows OS. All agents are simulated on a single machine.\footnote{Our implementation is available at \url{https://github.com/GlassBALDY/D-ripALM.git}.}

\paragraph{Restarted variant.}

The auxiliary variable $\vw^k$ in Algorithm \ref{alg:D-ripALM}
accumulates scaled subproblem errors through the update
\eqref{eq:D-ripALM-ywupdate}. Consequently, $\vw^k$ may drift
away from the candidate primal iterate $\vx^{k+1}$, potentially
increasing the cross term
$2|\langle \vw^k-\vx^{k+1},\sigma_k\Delta^{k+1}\rangle|$
in the relative-type error criterion \eqref{eq:D-ripALM-errcrit}.
A larger cross term may make the criterion more stringent,
requiring additional inner iterations to obtain an acceptable
subproblem solution.

Motivated by this observation, we consider a heuristic restarted
variant of D-ripALM, denoted by D-ripALM-re, following the
computational strategy of Eckstein and Silva
\cite[Section 6]{es2013practical}.
Similar to their implementation, we omit the cross term during
the initial outer iterations, but adopt a periodic restart
schedule rather than allowing only finitely many resets.
Specifically, we omit the cross term in
\eqref{eq:D-ripALM-errterm2} during the first three outer
iterations and restore it from the fourth iteration onward.
We then reset $\vw$ to the current primal iterate every two
outer iterations from the fourth through the tenth iteration,
and every three outer iterations thereafter. This restart schedule was selected based on preliminary numerical experiments to balance the additional inner iterations caused by a stringent error criterion and the frequency of auxiliary-variable resets.

Related restart mechanisms have also been studied for primal--dual methods, including the adaptive restarting scheme used in PDLP \cite{applegate2025pdlp}. Unlike these primal--dual restart schemes, our strategy resets only the auxiliary variable $\vw$ used in the relative-type error criterion, without restarting the primal--dual iterates or the momentum of the inner solver. Each agent performs this reset locally, without additional communication. Since omitting the cross term and periodically resetting $\vw$ modify the error criterion and the auxiliary-variable update in Algorithm \ref{alg:D-ripALM}, the convergence guarantees in Section \ref{sec:convergence} apply to the original D-ripALM algorithm but do not directly extend to D-ripALM-re.

\paragraph{Implementation details.}
For D-ripALM-re, we set $\vx^{0}=\vzero$ and $\vy^{0}=\vzero$. The local objectives of both test problems have the composite form $f_{i}=s_{i}+r_{i}$, where $s_{i}$ is smooth and the proximal mapping of $r_{i}$ can be evaluated efficiently. We solve each subproblem using FISTA \cite{beck2009fast}, initialized at the preceding outer iterate, with proximal gradient step size $(L+\sigma_{k}\lambda_{\max}(Z)+\tau_{k}/\sigma_{k})^{-1}$. Here, $L$ is the largest Lipschitz constant among the local gradients $\nabla s_{i}$. To evaluate the subproblem residual, we obtain a subgradient of $r_{i}$ at the candidate point from the optimality condition of the proximal step and add $\nabla s_{i}$ at that point. The resulting subgradient of $f_{i}$ is used in \eqref{eq:grad-natural} to compute $\Delta_{i}^{k+1}$. We check the inner stopping criterion every ten inner iterations and at the inner iteration limit.

\paragraph{Computational and communication costs.}
The main operations in an iteration of a single-loop method or an inner iteration of a double-loop method are a local gradient evaluation and a vector exchange with neighbors, with an additional proximal evaluation for LASSO. With sparse matrix--vector products, each update costs $\mathcal{O}(\operatorname{nnz}(A_{i})+m_{i}+(|\mathcal{N}_{i}|+1)d)$ operations per agent, excluding preprocessing and stopping tests. Here, $A_{i}$ is the local data matrix with $m_{i}$ rows and $\operatorname{nnz}(A_{i})$ nonzero entries. The same bound applies to storage per agent, excluding iteration histories. One exchange with all neighbors transmits $\mathcal{O}(|\mathcal{N}_{i}|d)$ scalars per agent, while stopping tests require additional residual evaluations and network aggregation. For the double-loop methods, the cost per outer iteration therefore depends on both the number of inner iterations and the frequency of inner stopping checks.

%%%%%%%%%%%%%%%%%%%%%%%%%%%%%%%%%%%%%%%%%%%%%%%%%%%%%%%%%%%%%%%%%%%%%%%%%
\subsection{Comparison with IDEAL}\label{sec:comp_diff_paras}

This subsection compares D-ripALM-re and IDEAL on distributed $\ell_2$-regularized logistic regression with synthetic data and MNIST.

We generate an Erd\H{o}s--R\'enyi graph on $n$ agents by including each edge independently with probability $0.2$, and resample the graph until it is connected. The mixing matrix $W$ is constructed using the Metropolis rule. Each agent $i$ has access to $m_i$ local training samples $\left(\bm{a}_{i j}, y_{i j}\right) \in \mathbb{R}^{d} \times\{-1,+1\}$. The resulting global optimization problem is
\begin{equation}\label{eq:log-reg}
\min\limits_{\bm{x} \in \mathbb{R}^d} \quad \sum_{i=1}^{n} \left[\sum_{j=1}^{m_i} \log \left(1+\exp \big(-y_{i j} \cdot \bm{a}_{i j}^{\mathsf{T}} \bm{x}\big)\right) + \frac{\lambda}{2}\|\bm{x}\|_2^2\right].
\end{equation}

We set $n=10$ and $\lambda=10^{-2}$ for both datasets described below:
\begin{itemize}
\item \textbf{Synthetic dataset.} We generate a global matrix $A=U\Sigma V^{\mathsf{T}}\in\RR^{m\times d}$, where $m=\sum_{i=1}^{n}m_{i}$. The orthogonal matrices $U\in\RR^{m\times m}$ and $V\in\RR^{d\times d}$ are obtained by QR factorization of matrices with independent standard Gaussian entries. The matrix $\Sigma\in\RR^{m\times d}$ has zero off-diagonal entries and diagonal entries evenly spaced from $\sqrt{L_{\max}}$ to zero, with $L_{\max}=10^{3}$, so that $\|A\|_{2}^{2}=L_{\max}$. We divide the rows of $A$ into $n$ consecutive blocks of equal size, one for each agent. The labels are generated as $y_{ij}=\operatorname{sign}(\bm{a}_{ij}^{\mathsf{T}}\bm{x}_{\mathrm{true}}+\epsilon_{ij})$, where $\bm{x}_{\mathrm{true}}=\bm{z}/\|\bm{z}\|$ with $\bm{z}\sim\mathcal{N}(0,I_{d})$, and the noise variables $\epsilon_{ij}$ are independent and follow $\mathcal{N}(0,0.1^{2})$. We set $d=1000$ and $m_i=40$.

\item \textbf{MNIST dataset.} We also consider binary classification of digits $0$ and $8$ using the MNIST training set. The resulting subset contains $m=11774$ samples, with labels $-1$ and $+1$ assigned to digits $0$ and $8$, respectively. Each $28\times28$ image is flattened into a feature vector of dimension $d=784$. Pixel values are divided by $255$, and each feature vector is then normalized to unit Euclidean norm. We conduct $10$ trials with different random seeds. In each trial, we randomly permute the sample-label pairs and partition them into $n$ consecutive blocks, one for each agent. The first $n-1$ blocks each contain $\lfloor m/n\rfloor$ samples, and the last block contains all remaining samples. We then generate a connected Erd\H{o}s--R\'enyi graph as described above.
\end{itemize}

IDEAL \cite[Algorithm 3]{arjevani2020ideal} employs an absolute-type error criterion to control the accuracy of solving the ALM subproblems. Since the local objectives are $\lambda$-strongly convex, the following gradient-residual test is sufficient to meet this criterion:
\begin{equation*}
\big\|\vx^{k+1} - \vx^{k+1,*}\big\|^2 \le \frac{1}{\lambda^{2}}\big\|\nabla\mathcal{L}_{\sigma_{k}}(\vx^{k+1},\vy^{k})\big\|^2 \leq \varepsilon_{k},
\end{equation*}
where $\vx^{k+1,*}$ denotes an exact solution of the ALM-subproblem at the $k$-th outer loop. The choice of the tolerance sequence $\{\varepsilon_k\}$ plays an important role: overly loose tolerances may hinder convergence, while overly stringent tolerances can significantly increase computational cost. Following the discussions in \cite{arjevani2020ideal}, we set $\varepsilon_k = \varepsilon_{0} \alpha^k$, with $\varepsilon_{0}\in\{10^{-2},10^{-1},10^{0},10^{1}\}$ and $\alpha\in\{0.2,0.4,0.6,0.8\}$, resulting in a total of $16$ parameter configurations. For IDEAL, we use its publicly available code\footnote{Code is available from the \href{https://proceedings.neurips.cc/paper/2020/hash/ed77eab0b8ff85d0a6a8365df1846978-Abstract.html}{IDEAL proceedings page}.}, retaining the accelerated gradient subsolver and its default parameter settings. We adapt the inner stopping test to the above error criterion and check it every ten inner iterations and at the inner iteration limit.

In contrast, D-ripALM-re adopts the relative-type error criterion \eqref{eq:D-ripALM-errcrit} with a restart strategy, using a single tolerance parameter $\rho \in [0,1)$. For comparisons, we consider $\rho \in \{0.01,\,0.05,\,0.1,\,0.3,\allowbreak\,0.5,\,0.7,\,0.9,\,0.99\}$, yielding eight choices. Moreover, we set the proximal parameter $\tau_k \equiv 10^{-3}$, adopt an increasing penalty parameter sequence $\sigma_k=\min\{1.5^{k}, 10^{4}\}$, and employ FISTA \cite{beck2009fast} as the subsolver.

Finally, we terminate all algorithms when the KKT residual at the current iterate $(\vx^{k+1},\vy^{k+1})$ satisfies
\begin{equation*}
\texttt{KKTres} \coloneqq \max\left\{\|\sqrt{Z}\vx^{k+1}\|, \,\big\|\nabla F(\vx^{k+1}) + \sqrt{Z}\vy^{k+1}\big\|\right\} < 10^{-6},
\end{equation*}
or when a maximum budget of 30,000 inner iterations is reached.

Table \ref{tab:log-reg} reports mean iteration counts, runtimes, and final KKT residuals over $10$ randomly generated synthetic instances. Here, \texttt{iter.(\#)} counts the total number of inner iterations for IDEAL and D-ripALM-re. These counts exclude communication for stopping tests and do not represent total communication rounds. The column \texttt{time(s)} reports runtime in seconds up to termination, including runs stopped at the iteration limit. Timings are measured in a single-machine simulation and do not include actual network delays.

With suitable tolerances, both methods require similar mean numbers of inner iterations. Across the tested settings, the smallest mean runtime is $1.44$ seconds for D-ripALM-re at $\rho=0.99$ and $0.54$ seconds for IDEAL at $(\varepsilon_{0},\alpha)=(0.1,0.2)$. At these settings, D-ripALM-re uses slightly fewer inner iterations on average, while IDEAL has a shorter mean runtime. D-ripALM-re uses a single tolerance parameter $\rho$, while the geometric tolerance sequence used by IDEAL has two parameters $\varepsilon_{0}$ and $\alpha$.

Table \ref{tab:log-reg-mnist} reports the mean MNIST results over $10$ trials. Within each trial, the same network and data partition are used for all tested parameter settings. IDEAL at $(\varepsilon_{0},\alpha)=(10,0.2)$ has a lower mean inner iteration count and a shorter mean runtime than D-ripALM-re at any tested setting. The lowest mean iteration count for D-ripALM-re is attained at $\rho=0.99$ on synthetic data and $\rho=0.9$ on MNIST, whereas the corresponding IDEAL setting changes from $(0.1,0.2)$ to $(10,0.2)$. The relative-type criterion requires only one tolerance parameter, and larger $\rho$ values generally yield fewer inner iterations and shorter runtimes in these tests.

\begin{table}[ht]
\centering
\caption{Numerical comparison for distributed $\ell_2$-regularized logistic regression on synthetic data. The total number of inner iterations, the outer iteration count in parentheses, and the total runtime \texttt{time(s)} are reported as mean $\pm$ sample standard deviation over $10$ instances. \texttt{KKTres} is the mean final KKT residual, reflecting both consensus and stationarity.}\label{tab:log-reg}
\begingroup
\small
\setlength{\tabcolsep}{5pt}
\renewcommand{\arraystretch}{1.05}
\resizebox{\linewidth}{!}{%
\shortstack[c]{%
\begin{tabular}{@{}llll@{\hspace{12pt}}llll@{}}
\toprule
\multicolumn{8}{c}{IDEAL} \\
\cmidrule{1-8}
$(\varepsilon_0,\alpha)$ & \texttt{iter.(\#)} & \texttt{time(s)} & \texttt{KKTres} & $(\varepsilon_0,\alpha)$ & \texttt{iter.(\#)} & \texttt{time(s)} & \texttt{KKTres} \\
\cmidrule(lr){1-4}\cmidrule(lr){5-8}
(0.01, 0.2) & 3895 $\pm$ 134 (9 $\pm$ 0) & 0.57 $\pm$ 0.04 & 7.0e-07 & (1, 0.2) & 4019 $\pm$ 37 (12 $\pm$ 0) & 0.58 $\pm$ 0.01 & 6.3e-07 \\
(0.01, 0.4) & 4054 $\pm$ 76 (16 $\pm$ 0) & 0.59 $\pm$ 0.02 & 7.5e-07 & (1, 0.4) & 4687 $\pm$ 44 (21 $\pm$ 0) & 0.70 $\pm$ 0.03 & 6.5e-07 \\
(0.01, 0.6) & 4460 $\pm$ 75 (27 $\pm$ 0) & 0.67 $\pm$ 0.03 & 9.4e-07 & (1, 0.6) & 5459 $\pm$ 124 (37 $\pm$ 1) & 0.81 $\pm$ 0.03 & 8.6e-07 \\
(0.01, 0.8) & 5556 $\pm$ 64 (62 $\pm$ 0) & 0.86 $\pm$ 0.04 & 9.8e-07 & (1, 0.8) & 7222 $\pm$ 87 (83 $\pm$ 0) & 1.11 $\pm$ 0.03 & 9.4e-07 \\
(0.1, 0.2) & 3740 $\pm$ 157 (10 $\pm$ 1) & 0.54 $\pm$ 0.04 & 7.1e-07 & (10, 0.2) & 4125 $\pm$ 42 (13 $\pm$ 0) & 0.60 $\pm$ 0.02 & 8.8e-07 \\
(0.1, 0.4) & 4196 $\pm$ 44 (18 $\pm$ 0) & 0.63 $\pm$ 0.03 & 8.2e-07 & (10, 0.4) & 4907 $\pm$ 52 (23 $\pm$ 0) & 0.73 $\pm$ 0.03 & 8.3e-07 \\
(0.1, 0.6) & 5021 $\pm$ 95 (32 $\pm$ 0) & 0.75 $\pm$ 0.04 & 8.8e-07 & (10, 0.6) & 5908 $\pm$ 58 (41 $\pm$ 0) & 0.89 $\pm$ 0.03 & 8.8e-07 \\
(0.1, 0.8) & 6647 $\pm$ 79 (73 $\pm$ 0) & 1.03 $\pm$ 0.04 & 9.1e-07 & (10, 0.8) & 7921 $\pm$ 79 (93 $\pm$ 0) & 1.23 $\pm$ 0.04 & 9.7e-07 \\
\bottomrule
\end{tabular} \\[8pt]
\begin{tabular}{@{}llll@{}}
\toprule
\multicolumn{4}{c}{D-ripALM-re} \\
\cmidrule{1-4}
$\rho$ & \texttt{iter.(\#)} & \texttt{time(s)} & \texttt{KKTres} \\
\midrule
0.99 & 3478 $\pm$ 682 (11 $\pm$ 0) & 1.44 $\pm$ 0.31 & 4.9e-07 \\
0.9 & 3779 $\pm$ 1049 (11 $\pm$ 0) & 1.55 $\pm$ 0.43 & 4.3e-07 \\
0.7 & 4113 $\pm$ 630 (11 $\pm$ 0) & 1.68 $\pm$ 0.23 & 1.7e-07 \\
0.5 & 4097 $\pm$ 290 (11 $\pm$ 0) & 1.67 $\pm$ 0.12 & 2.4e-07 \\
0.3 & 4250 $\pm$ 387 (10 $\pm$ 1) & 1.74 $\pm$ 0.16 & 4.4e-07 \\
0.1 & 5375 $\pm$ 319 (10 $\pm$ 0) & 2.19 $\pm$ 0.13 & 2.0e-07 \\
0.05 & 6262 $\pm$ 964 (10 $\pm$ 0) & 2.57 $\pm$ 0.46 & 1.7e-07 \\
0.01 & 7900 $\pm$ 754 (10 $\pm$ 1) & 3.24 $\pm$ 0.30 & 4.0e-07 \\
\bottomrule
\end{tabular}
}}
\endgroup
\end{table}

\begin{table}[ht]
\centering
\caption{Numerical comparison for distributed $\ell_2$-regularized logistic regression on MNIST. The total number of inner iterations, the outer iteration count in parentheses, and the total runtime \texttt{time(s)} are reported as mean $\pm$ sample standard deviation over $10$ trials. \texttt{KKTres} is the mean final KKT residual, reflecting both consensus and stationarity.}\label{tab:log-reg-mnist}
\begingroup
\small
\setlength{\tabcolsep}{5pt}
\renewcommand{\arraystretch}{1.05}
\resizebox{\linewidth}{!}{%
\shortstack[c]{%
\begin{tabular}{@{}llll@{\hspace{12pt}}llll@{}}
\toprule
\multicolumn{8}{c}{IDEAL} \\
\cmidrule{1-8}
$(\varepsilon_0,\alpha)$ & \texttt{iter.(\#)} & \texttt{time(s)} & \texttt{KKTres} & $(\varepsilon_0,\alpha)$ & \texttt{iter.(\#)} & \texttt{time(s)} & \texttt{KKTres} \\
\cmidrule(lr){1-4}\cmidrule(lr){5-8}
(0.01, 0.2) & 6339 $\pm$ 1756 (10 $\pm$ 1) & 16.81 $\pm$ 4.64 & 5.3e-07 & (1, 0.2) & 4529 $\pm$ 479 (12 $\pm$ 1) & 12.18 $\pm$ 1.22 & 5.9e-07 \\
(0.01, 0.4) & 5740 $\pm$ 274 (16 $\pm$ 0) & 15.54 $\pm$ 0.70 & 6.5e-07 & (1, 0.4) & 5150 $\pm$ 45 (21 $\pm$ 0) & 14.18 $\pm$ 0.94 & 6.5e-07 \\
(0.01, 0.6) & 6209 $\pm$ 207 (28 $\pm$ 1) & 16.91 $\pm$ 0.57 & 8.6e-07 & (1, 0.6) & 5992 $\pm$ 98 (36 $\pm$ 1) & 16.30 $\pm$ 0.26 & 8.8e-07 \\
(0.01, 0.8) & 7497 $\pm$ 243 (62 $\pm$ 0) & 20.78 $\pm$ 0.70 & 9.8e-07 & (1, 0.8) & 7814 $\pm$ 71 (83 $\pm$ 0) & 21.88 $\pm$ 0.37 & 9.4e-07 \\
(0.1, 0.2) & 5219 $\pm$ 1113 (11 $\pm$ 1) & 13.90 $\pm$ 2.90 & 7.2e-07 & (10, 0.2) & 4468 $\pm$ 97 (13 $\pm$ 0) & 12.06 $\pm$ 0.39 & 8.4e-07 \\
(0.1, 0.4) & 5138 $\pm$ 234 (18 $\pm$ 0) & 13.79 $\pm$ 0.61 & 8.1e-07 & (10, 0.4) & 5234 $\pm$ 53 (23 $\pm$ 0) & 14.52 $\pm$ 0.89 & 8.2e-07 \\
(0.1, 0.6) & 5679 $\pm$ 265 (32 $\pm$ 0) & 15.45 $\pm$ 0.73 & 8.8e-07 & (10, 0.6) & 6212 $\pm$ 78 (41 $\pm$ 0) & 16.99 $\pm$ 0.23 & 8.8e-07 \\
(0.1, 0.8) & 7144 $\pm$ 218 (73 $\pm$ 0) & 19.92 $\pm$ 0.55 & 9.1e-07 & (10, 0.8) & 8212 $\pm$ 65 (93 $\pm$ 0) & 23.10 $\pm$ 0.26 & 9.8e-07 \\
\bottomrule
\end{tabular} \\[8pt]
\begin{tabular}{@{}llll@{}}
\toprule
\multicolumn{4}{c}{D-ripALM-re} \\
\cmidrule{1-4}
$\rho$ & \texttt{iter.(\#)} & \texttt{time(s)} & \texttt{KKTres} \\
\midrule
0.99 & 4714 $\pm$ 602 (15 $\pm$ 1) & 13.45 $\pm$ 1.66 & 5.7e-07 \\
0.9 & 4625 $\pm$ 481 (15 $\pm$ 1) & 13.22 $\pm$ 1.27 & 6.2e-07 \\
0.7 & 4774 $\pm$ 590 (15 $\pm$ 1) & 13.65 $\pm$ 1.56 & 5.9e-07 \\
0.5 & 4711 $\pm$ 428 (16 $\pm$ 1) & 13.45 $\pm$ 1.18 & 5.8e-07 \\
0.3 & 5135 $\pm$ 577 (16 $\pm$ 1) & 14.80 $\pm$ 1.51 & 5.9e-07 \\
0.1 & 6181 $\pm$ 1064 (16 $\pm$ 1) & 17.73 $\pm$ 2.93 & 5.2e-07 \\
0.05 & 6907 $\pm$ 840 (16 $\pm$ 1) & 19.85 $\pm$ 2.28 & 5.3e-07 \\
0.01 & 9640 $\pm$ 1533 (16 $\pm$ 1) & 27.68 $\pm$ 4.31 & 5.2e-07 \\
\bottomrule
\end{tabular}
}}
\endgroup
\end{table}

%%%%%%%%%%%%%%%%%%%%%%%%%%%%%%%%%%%%%%%%%%%%%%%%%%%%%%%%%%%%%%%%%%%%%%%%%
\subsection{Comparison with PG-EXTRA and NIDS}\label{sec:comp_diff_algos}

In this subsection, we compare D-ripALM-re with two representative single-loop distributed algorithms, PG-EXTRA and NIDS, on the distributed $\ell_1$-regularized least-squares regression (LASSO) problem.

We consider a connected network with $n$ agents. Each agent $i$ holds a feature matrix $A_{i}\in\RR^{m_i\times d}$ and a response vector $\vb_{i}\in\RR^{m_{i}}$. The global optimization problem is
\begin{equation}\label{eq:lasso}
\min\limits_{\bm{x} \in \mathbb{R}^d} \quad \sum_{i=1}^{n} \left[\frac{1}{2}\big\|A_{i}\bm{x} - \vb_{i}\big\|^{2} + \frac{\lambda}{n}\|\bm{x}\|_1\right].
\end{equation}
We use synthetic data and two real datasets from the LIBSVM collection.
\begin{itemize}
\item \textbf{Synthetic dataset.} We generate the global matrix $A$ as in Section \ref{sec:comp_diff_paras}, with $L_{\max}=10^{2}$, and partition its rows into $n$ consecutive blocks of equal size to obtain the local matrices $A_{i}$. The associated measurement $\vb_{i} \in \RR^{m_{i}}$ is constructed as $\vb_{i} = A_{i}\bm{x}_{\text{true}} + \bm{\epsilon}_{i}$, where $\bm{x}_{\text{true}}=\bm{z}/\|\bm{z}\|$ with $\bm{z}=\texttt{sprandn(d,1,0.1)}$ and $\bm{\epsilon}_{i} \sim \mathcal{N}(\vzero, 0.1^2 I_{m_{i}})$. In these experiments, we set $d=1000$ and $\sum_{i=1}^{n}m_{i} = 200$.

\item \textbf{Real datasets.} We test \texttt{abalone\_scale} ($m=4177$, $d=8$) and \texttt{space\_ga\_scale} ($m=3107$, $d=6$) from the LIBSVM collection.\footnote{Available at \url{https://www.csie.ntu.edu.tw/~cjlin/libsvmtools/datasets/}.} For each dataset and network topology, we conduct $10$ trials using the same shuffling and partitioning procedure as for MNIST in Section \ref{sec:comp_diff_paras}.
\end{itemize}

In the following experiments, we set $n=20$. The regularization parameter is chosen as $\lambda = \lambda_{c}\|A^{\mathsf{T}}\bm{b}\|_{\infty}$, where $A = [A_{1}; A_{2};\ldots;A_{n}]\in\RR^{m\times d}$ and $\bm{b} = [\vb_{1}; \vb_{2};\ldots;\vb_{n}]\in\RR^{m}$ are the stacked data, and $\lambda_{c}\in\{10^{-1}, 10^{-2}, 10^{-3}\}$. To examine the impact of network connectivity, we consider three commonly used communication topologies: a ring graph, an Erd\H{o}s--R\'enyi graph with connection probability 0.2, and a geometric random graph. For the geometric graph, we sample node positions uniformly in $[0,1]^{2}$ and connect two nodes if their Euclidean distance is less than $\sqrt{6/(\pi n)}$. We resample random graphs until they are connected. For the real datasets, the Erd\H{o}s--R\'enyi and geometric graphs are generated anew in each trial. For each topology, the mixing matrix $W$ is constructed using the Metropolis rule.

PG-EXTRA and NIDS are single-loop methods. We use their publicly available implementations\footnote{Available at \url{https://github.com/mingyan08/NIDS}.}, with the necessary adaptations for the test problems and stopping tests. We follow the parameter recommendations in \cite{shi2015proximal,li2019nids}.
For NIDS, we set $\alpha_{i}=1/L_{i}$ and $c=1/(2\max_{i}\alpha_{i})$, where $L_{i}=\|A_{i}\|_{2}^{2}$. For PG-EXTRA, we use the common step size $\alpha=0.99(1+\lambda_{\min}(W))/L$, with $L=\max_{i}L_{i}$. These parameters are computed from the smoothness constants and network information for each instance, without an empirical parameter search. For D-ripALM-re, we use the settings from Section \ref{sec:comp_diff_paras} with $\rho=0.99$. All algorithms are terminated either when the maximum number of 30,000 (inner) iterations is reached or when the current iterate $\vx^{k+1}$ satisfies
\begin{equation*}
\texttt{KKTres} \coloneqq \max\left\{\big\|\sqrt{Z}\vx^{k+1}\big\|, \,\frac{\big\|\bar{\bm{x}}^{k+1} - \mathtt{prox}_{\lambda\|\cdot\|_{1}}\big(\bar{\bm{x}}^{k+1} - A^{\mathsf{T}}(A\bar{\bm{x}}^{k+1}-\bm{b})\big)\big\|}{1 + \|A\bar{\bm{x}}^{k+1}-\bm{b}\| + \|\bar{\bm{x}}^{k+1}\|}\right\} < 10^{-6},
\end{equation*}
where $\bar{\bm{x}}^{k+1} = \frac{1}{n}\big(\vone_{n}^{\mathsf{T}}\otimes I_{d}\big)\vx^{k+1}$ is the averaged solution over all agents. Since the averaged solution $\bar{\bm{x}}^{k+1}$ is typically dense, we additionally perform a simple refinement step by setting $\bar{x}_{i}^{k+1} = 0$ if $|\bar{x}_{i}^{k+1}| / \|\bar{\bm{x}}^{k+1}\|_{\infty} < 10^{-8}$ to improve the quality of the approximate solution.

Table \ref{tab:LASSO} reports mean iteration counts, runtimes, and final KKT residuals over $10$ randomly generated synthetic instances. Here, \texttt{iter.(\#)} counts inner iterations for D-ripALM-re and iterations for PG-EXTRA and NIDS. For $\lambda_{c}=10^{-1}$ and $10^{-2}$, D-ripALM-re has lower mean iteration counts and shorter mean runtimes than PG-EXTRA and NIDS on all three topologies. Its mean final residuals are below $10^{-6}$, whereas those of the two single-loop methods remain above $10^{-6}$. At $\lambda_{c}=10^{-3}$, all three methods reach the iteration limit, but D-ripALM-re attains much smaller mean final KKT residuals with comparable runtimes.

Tables \ref{tab:lasso_libsvm_abalone_scale} and \ref{tab:lasso_libsvm_space_ga_scale} summarize the results over $10$ trials on the two real datasets. On \texttt{abalone\_scale}, D-ripALM-re has lower mean iteration counts and shorter mean runtimes than both PG-EXTRA and NIDS in all tested settings. In particular, at $\lambda_{c}=10^{-3}$, it attains mean final KKT residuals below $10^{-6}$ on all three topologies, while both single-loop methods reach the iteration limit with mean final residuals above $10^{-6}$. On \texttt{space\_ga\_scale}, D-ripALM-re also has lower mean iteration counts and shorter mean runtimes at $\lambda_{c}=10^{-2}$ and $10^{-3}$ on all three topologies. At $\lambda_{c}=10^{-1}$, PG-EXTRA has a shorter mean runtime on the ring and geometric graphs, whereas the three methods have very similar mean runtimes on the Erd\H{o}s--R\'enyi graph ($0.15$ seconds for D-ripALM-re versus $0.14$ seconds for NIDS and PG-EXTRA). These results suggest that D-ripALM-re can offer competitive computational efficiency and solution accuracy on the tested real datasets.

\begin{table}[htbp]
\centering
\caption{Numerical comparison for distributed LASSO on synthetic data. The iteration counts, the outer iteration count of D-ripALM-re in parentheses, and the total runtime \texttt{time(s)} are reported as mean $\pm$ sample standard deviation over $10$ instances. \texttt{KKTres} is the mean final KKT residual, reflecting both consensus and stationarity.}\label{tab:LASSO}
\begingroup
\small
\setlength{\tabcolsep}{5pt}
\renewcommand{\arraystretch}{1.05}
\resizebox{\linewidth}{!}{%
\begin{tabular}{llllll}
\toprule
$\lambda_{c}$ & Topology & Method & \texttt{iter.(\#)} & \texttt{time(s)} & \texttt{KKTres} \\
\midrule
\multirow{9}{*}[-5.4pt]{$10^{-1}$} & \multirow{3}{*}{Ring} & D-ripALM-re & 8601 $\pm$ 1134 (20 $\pm$ 0) & 7.00 $\pm$ 0.90 & 5.0e-07 \\
&  & NIDS & 28995 $\pm$ 1979 & 28.23 $\pm$ 2.05 & 7.4e-06 \\
&  & PG-EXTRA & 30000 $\pm$ 0 & 29.22 $\pm$ 0.87 & 9.2e-05 \\
\cmidrule(lr){2-6}
& \multirow{3}{*}{Erd\H{o}s--R\'enyi} & D-ripALM-re & 6665 $\pm$ 1890 (19 $\pm$ 1) & 6.05 $\pm$ 1.61 & 4.7e-07 \\
&  & NIDS & 29007 $\pm$ 1954 & 28.14 $\pm$ 1.79 & 7.4e-06 \\
&  & PG-EXTRA & 30000 $\pm$ 0 & 29.18 $\pm$ 0.62 & 6.3e-05 \\
\cmidrule(lr){2-6}
& \multirow{3}{*}{Geometric} & D-ripALM-re & 8480 $\pm$ 1934 (20 $\pm$ 0) & 7.91 $\pm$ 1.57 & 3.4e-07 \\
&  & NIDS & 28995 $\pm$ 1981 & 28.04 $\pm$ 1.91 & 7.4e-06 \\
&  & PG-EXTRA & 30000 $\pm$ 0 & 28.94 $\pm$ 0.30 & 2.9e-05 \\
\midrule
\multirow{9}{*}[-5.4pt]{$10^{-2}$} & \multirow{3}{*}{Ring} & D-ripALM-re & 29825 $\pm$ 484 (20 $\pm$ 0) & 24.00 $\pm$ 0.43 & 6.7e-07 \\
&  & NIDS & 30000 $\pm$ 0 & 28.84 $\pm$ 0.24 & 2.5e-03 \\
&  & PG-EXTRA & 30000 $\pm$ 0 & 28.75 $\pm$ 0.21 & 4.7e-03 \\
\cmidrule(lr){2-6}
& \multirow{3}{*}{Erd\H{o}s--R\'enyi} & D-ripALM-re & 24910 $\pm$ 4387 (18 $\pm$ 1) & 22.53 $\pm$ 4.27 & 2.7e-07 \\
&  & NIDS & 30000 $\pm$ 0 & 28.75 $\pm$ 0.19 & 2.4e-03 \\
&  & PG-EXTRA & 30000 $\pm$ 0 & 28.79 $\pm$ 0.19 & 4.3e-03 \\
\cmidrule(lr){2-6}
& \multirow{3}{*}{Geometric} & D-ripALM-re & 26842 $\pm$ 4242 (19 $\pm$ 0) & 24.84 $\pm$ 3.95 & 4.1e-07 \\
&  & NIDS & 30000 $\pm$ 0 & 28.83 $\pm$ 0.20 & 2.4e-03 \\
&  & PG-EXTRA & 30000 $\pm$ 0 & 28.75 $\pm$ 0.25 & 3.5e-03 \\
\midrule
\multirow{9}{*}[-5.4pt]{$10^{-3}$} & \multirow{3}{*}{Ring} & D-ripALM-re & 30000 $\pm$ 0 (17 $\pm$ 1) & 24.17 $\pm$ 0.18 & 2.0e-04 \\
&  & NIDS & 30000 $\pm$ 0 & 28.81 $\pm$ 0.20 & 5.0e-03 \\
&  & PG-EXTRA & 30000 $\pm$ 0 & 28.70 $\pm$ 0.18 & 7.4e-03 \\
\cmidrule(lr){2-6}
& \multirow{3}{*}{Erd\H{o}s--R\'enyi} & D-ripALM-re & 30000 $\pm$ 0 (16 $\pm$ 0) & 27.08 $\pm$ 0.83 & 3.8e-05 \\
&  & NIDS & 30000 $\pm$ 0 & 28.97 $\pm$ 0.31 & 5.0e-03 \\
&  & PG-EXTRA & 30000 $\pm$ 0 & 28.85 $\pm$ 0.30 & 6.9e-03 \\
\cmidrule(lr){2-6}
& \multirow{3}{*}{Geometric} & D-ripALM-re & 30000 $\pm$ 0 (16 $\pm$ 0) & 27.76 $\pm$ 0.69 & 1.2e-04 \\
&  & NIDS & 30000 $\pm$ 0 & 28.81 $\pm$ 0.19 & 5.0e-03 \\
&  & PG-EXTRA & 30000 $\pm$ 0 & 28.78 $\pm$ 0.31 & 6.2e-03 \\
\bottomrule
\end{tabular}}
\endgroup
\end{table}

\begin{table}[htbp]
\centering
\caption{Numerical comparison for distributed LASSO on \texttt{abalone\_scale} ($m=4177$, $d=8$) over three network topologies. The iteration counts, the outer iteration count of D-ripALM-re in parentheses, and the total runtime \texttt{time(s)} are reported as mean $\pm$ sample standard deviation over $10$ trials. \texttt{KKTres} is the mean final KKT residual, reflecting both consensus and stationarity.}\label{tab:lasso_libsvm_abalone_scale}
\begingroup
\small
\setlength{\tabcolsep}{5pt}
\renewcommand{\arraystretch}{1.05}
\resizebox{\linewidth}{!}{%
\begin{tabular}{lllllc}
\toprule
$\lambda_c$ & Topology & Method & \texttt{iter.(\#)} & \texttt{time(s)} & \texttt{KKTres} \\
\midrule
\multirow{9}{*}[-5.4pt]{$10^{-1}$} & \multirow{3}{*}{Ring} & D-ripALM-re & 7006 $\pm$ 908 (31 $\pm$ 1) & 0.87 $\pm$ 0.11 & 5.2e-07 \\
&  & NIDS & 13764 $\pm$ 16 & 5.35 $\pm$ 0.08 & 1.0e-06 \\
&  & PG-EXTRA & 22951 $\pm$ 378 & 5.07 $\pm$ 0.16 & 1.0e-06 \\
\cmidrule(lr){2-6}
& \multirow{3}{*}{Erd\H{o}s--R\'enyi} & D-ripALM-re & 3809 $\pm$ 687 (25 $\pm$ 1) & 0.56 $\pm$ 0.09 & 4.2e-07 \\
&  & NIDS & 13764 $\pm$ 16 & 5.36 $\pm$ 0.08 & 1.0e-06 \\
&  & PG-EXTRA & 20886 $\pm$ 811 & 4.60 $\pm$ 0.21 & 1.0e-06 \\
\cmidrule(lr){2-6}
& \multirow{3}{*}{Geometric} & D-ripALM-re & 5716 $\pm$ 926 (29 $\pm$ 1) & 0.88 $\pm$ 0.14 & 5.2e-07 \\
&  & NIDS & 13765 $\pm$ 16 & 5.37 $\pm$ 0.09 & 1.0e-06 \\
&  & PG-EXTRA & 17900 $\pm$ 410 & 3.93 $\pm$ 0.11 & 1.0e-06 \\
\midrule
\multirow{9}{*}[-5.4pt]{$10^{-2}$} & \multirow{3}{*}{Ring} & D-ripALM-re & 5190 $\pm$ 394 (34 $\pm$ 1) & 0.65 $\pm$ 0.05 & 7.2e-07 \\
&  & NIDS & 5581 $\pm$ 7 & 2.20 $\pm$ 0.04 & 1.0e-06 \\
&  & PG-EXTRA & 9310 $\pm$ 153 & 2.06 $\pm$ 0.06 & 1.0e-06 \\
\cmidrule(lr){2-6}
& \multirow{3}{*}{Erd\H{o}s--R\'enyi} & D-ripALM-re & 3318 $\pm$ 464 (26 $\pm$ 1) & 0.50 $\pm$ 0.06 & 4.6e-07 \\
&  & NIDS & 5582 $\pm$ 7 & 2.17 $\pm$ 0.03 & 1.0e-06 \\
&  & PG-EXTRA & 8472 $\pm$ 329 & 1.89 $\pm$ 0.08 & 1.0e-06 \\
\cmidrule(lr){2-6}
& \multirow{3}{*}{Geometric} & D-ripALM-re & 4455 $\pm$ 461 (31 $\pm$ 2) & 0.69 $\pm$ 0.07 & 6.3e-07 \\
&  & NIDS & 5581 $\pm$ 6 & 2.18 $\pm$ 0.05 & 1.0e-06 \\
&  & PG-EXTRA & 7260 $\pm$ 167 & 1.60 $\pm$ 0.04 & 1.0e-06 \\
\midrule
\multirow{9}{*}[-5.4pt]{$10^{-3}$} & \multirow{3}{*}{Ring} & D-ripALM-re & 12223 $\pm$ 1228 (37 $\pm$ 1) & 1.52 $\pm$ 0.15 & 7.9e-07 \\
&  & NIDS & 30000 $\pm$ 0 & 11.66 $\pm$ 0.14 & 1.1e-06 \\
&  & PG-EXTRA & 30000 $\pm$ 0 & 6.58 $\pm$ 0.09 & 1.5e-04 \\
\cmidrule(lr){2-6}
& \multirow{3}{*}{Erd\H{o}s--R\'enyi} & D-ripALM-re & 7121 $\pm$ 970 (27 $\pm$ 1) & 1.05 $\pm$ 0.14 & 6.1e-07 \\
&  & NIDS & 30000 $\pm$ 0 & 11.61 $\pm$ 0.13 & 1.1e-06 \\
&  & PG-EXTRA & 30000 $\pm$ 0 & 6.61 $\pm$ 0.09 & 7.4e-05 \\
\cmidrule(lr){2-6}
& \multirow{3}{*}{Geometric} & D-ripALM-re & 9786 $\pm$ 1596 (34 $\pm$ 2) & 1.49 $\pm$ 0.23 & 7.9e-07 \\
&  & NIDS & 30000 $\pm$ 0 & 11.61 $\pm$ 0.15 & 1.1e-06 \\
&  & PG-EXTRA & 30000 $\pm$ 0 & 6.63 $\pm$ 0.13 & 1.9e-05 \\
\bottomrule
\end{tabular}}
\endgroup
\end{table}

\begin{table}[htbp]
\centering
\caption{Numerical comparison for distributed LASSO on \texttt{space\_ga\_scale} ($m=3107$, $d=6$), with the same reporting convention as in Table \ref{tab:lasso_libsvm_abalone_scale}.}\label{tab:lasso_libsvm_space_ga_scale}
\begingroup
\small
\setlength{\tabcolsep}{5pt}
\renewcommand{\arraystretch}{1.05}
\resizebox{\linewidth}{!}{%
\begin{tabular}{lllllc}
\toprule
$\lambda_c$ & Topology & Method & \texttt{iter.(\#)} & \texttt{time(s)} & \texttt{KKTres} \\
\midrule
\multirow{9}{*}[-5.4pt]{$10^{-1}$} & \multirow{3}{*}{Ring} & D-ripALM-re & 1967 $\pm$ 338 (23 $\pm$ 1) & 0.24 $\pm$ 0.04 & 5.4e-07 \\
&  & NIDS & 550 $\pm$ 18 & 0.20 $\pm$ 0.01 & 9.9e-07 \\
&  & PG-EXTRA & 743 $\pm$ 43 & 0.15 $\pm$ 0.01 & 9.9e-07 \\
\cmidrule(lr){2-6}
& \multirow{3}{*}{Erd\H{o}s--R\'enyi} & D-ripALM-re & 992 $\pm$ 226 (20 $\pm$ 1) & 0.15 $\pm$ 0.03 & 3.3e-07 \\
&  & NIDS & 383 $\pm$ 2 & 0.14 $\pm$ 0.00 & 9.7e-07 \\
&  & PG-EXTRA & 677 $\pm$ 56 & 0.14 $\pm$ 0.01 & 9.9e-07 \\
\cmidrule(lr){2-6}
& \multirow{3}{*}{Geometric} & D-ripALM-re & 1402 $\pm$ 374 (22 $\pm$ 1) & 0.22 $\pm$ 0.06 & 4.2e-07 \\
&  & NIDS & 441 $\pm$ 81 & 0.16 $\pm$ 0.03 & 9.9e-07 \\
&  & PG-EXTRA & 578 $\pm$ 34 & 0.11 $\pm$ 0.01 & 9.9e-07 \\
\midrule
\multirow{9}{*}[-5.4pt]{$10^{-2}$} & \multirow{3}{*}{Ring} & D-ripALM-re & 8442 $\pm$ 1069 (24 $\pm$ 0) & 1.02 $\pm$ 0.13 & 3.9e-07 \\
&  & NIDS & 12974 $\pm$ 73 & 4.63 $\pm$ 0.06 & 1.0e-06 \\
&  & PG-EXTRA & 24867 $\pm$ 1425 & 4.88 $\pm$ 0.29 & 1.0e-06 \\
\cmidrule(lr){2-6}
& \multirow{3}{*}{Erd\H{o}s--R\'enyi} & D-ripALM-re & 4119 $\pm$ 1138 (21 $\pm$ 1) & 0.60 $\pm$ 0.16 & 5.7e-07 \\
&  & NIDS & 12975 $\pm$ 73 & 4.60 $\pm$ 0.06 & 1.0e-06 \\
&  & PG-EXTRA & 22654 $\pm$ 1864 & 4.44 $\pm$ 0.38 & 1.0e-06 \\
\cmidrule(lr){2-6}
& \multirow{3}{*}{Geometric} & D-ripALM-re & 5827 $\pm$ 1576 (23 $\pm$ 1) & 0.87 $\pm$ 0.22 & 4.2e-07 \\
&  & NIDS & 12974 $\pm$ 73 & 4.61 $\pm$ 0.07 & 1.0e-06 \\
&  & PG-EXTRA & 19392 $\pm$ 1107 & 3.82 $\pm$ 0.24 & 1.0e-06 \\
\midrule
\multirow{9}{*}[-5.4pt]{$10^{-3}$} & \multirow{3}{*}{Ring} & D-ripALM-re & 7953 $\pm$ 876 (24 $\pm$ 0) & 0.97 $\pm$ 0.11 & 3.7e-07 \\
&  & NIDS & 11926 $\pm$ 66 & 4.26 $\pm$ 0.07 & 1.0e-06 \\
&  & PG-EXTRA & 22859 $\pm$ 1310 & 4.52 $\pm$ 0.31 & 1.0e-06 \\
\cmidrule(lr){2-6}
& \multirow{3}{*}{Erd\H{o}s--R\'enyi} & D-ripALM-re & 4618 $\pm$ 1235 (21 $\pm$ 1) & 0.67 $\pm$ 0.18 & 3.6e-07 \\
&  & NIDS & 11928 $\pm$ 67 & 4.27 $\pm$ 0.07 & 1.0e-06 \\
&  & PG-EXTRA & 20826 $\pm$ 1713 & 4.10 $\pm$ 0.35 & 1.0e-06 \\
\cmidrule(lr){2-6}
& \multirow{3}{*}{Geometric} & D-ripALM-re & 5994 $\pm$ 1845 (23 $\pm$ 1) & 0.90 $\pm$ 0.27 & 4.4e-07 \\
&  & NIDS & 11926 $\pm$ 67 & 4.26 $\pm$ 0.08 & 1.0e-06 \\
&  & PG-EXTRA & 17826 $\pm$ 1018 & 3.50 $\pm$ 0.22 & 1.0e-06 \\
\bottomrule
\end{tabular}}
\endgroup
\end{table}

%%%%%%%%%%%%%%%%%%%%%%%%%%%%%%%%%%%%%%%%%%%%%%%%%%%%%%%%%%%%%%%%%%%%%%%%%
\section{Conclusions}\label{sec:conclusion}

This paper proposes D-ripALM, a double-loop distributed ALM-type algorithm for decentralized consensus optimization. D-ripALM (Algorithm \ref{alg:D-ripALM}) builds on its centralized counterpart ripALM \cite{yang2025ripalm,zhu2024ripalm} and inherits rigorous convergence guarantees under the stated assumptions. Its double-loop structure permits efficient, problem-tailored subsolvers for both smooth and nonsmooth convex problems in distributed settings. The relative-type error criterion provides an adaptive mechanism for switching between (approximate) subproblem solves and multiplier updates, using one relative tolerance parameter instead of a prescribed absolute tolerance sequence. The experiments on synthetic and real datasets suggest that the restarted variant D-ripALM-re can offer competitive performance relative to representative decentralized algorithms. Future work may explore applications to time-varying and/or directed networks, asynchronous communication, and extensions to distributed resource allocation problems.

%%%%%%%%%%%%%%%%%%%%%%%%%%%%%%%%%%%%%%%%%%%%%%%%%%%%%%%%%%%%%%%%%%%%%%%%%
\subsection*{Acknowledgment}

The authors thank the editor and anonymous reviewers for their valuable comments and suggestions, which have helped to improve the paper. The research of Hong Wang is supported in part by National Natural Science Foundation of China under grant 12201428. The research of Lei Yang is supported in part by the National Natural Science Foundation of China under grant 12301411, and the Natural Science Foundation of Guangdong under grant 2026A1515010350.

\section*{Declarations}

\noindent\textbf{Conflict of interest}\quad
The authors declare that they have no conflict of interest.

\medskip
\noindent\textbf{Data availability}\quad
The synthetic test instances were generated as described in Sections \ref{sec:comp_diff_paras} and \ref{sec:comp_diff_algos}. The MNIST dataset is publicly available from the CVDF mirror at \url{https://github.com/cvdfoundation/mnist}. The \texttt{abalone\_scale} and \texttt{space\_ga\_scale} datasets are publicly available from the LIBSVM collection at \url{https://www.csie.ntu.edu.tw/~cjlin/libsvmtools/datasets/regression.html}.

%%%%%%%%%%%%%%%%%%%%%%%%%%%%%%%%%%%%%%%%%%%%%%%%%%%%%%%%%%%%%%%%%%%%%%%%%
\clearpage
\bibliographystyle{plain}
\bibliography{Ref_d_ripALM}

\end{document}